\documentstyle[12pt,amssymb,amscd]{amsart}  % amslatex

\theoremstyle{remark}

\theoremstyle{definition}

\numberwithin{equation}{subsection}

\newcommand{\bbC}{{\Bbb C}}

\newcommand{\bbZ}{{\Bbb Z}}

\newcommand{\cO}{{\cal O}}

\newcommand{\cV}{{\cal V}}

\newcommand{\cI}{{\cal I}}
\newcommand{\cC}{{\cal C}}
\newcommand{\cE}{{\cal E}}

\newcommand{\cU}{{\cal U}}

\newcommand{\cT}{{\cal T}}
\newcommand{\cQ}{{\cal Q}}

\newcommand{\Hom}{\operatorname{Hom}}

\newcommand{\isomo}{\overset{\sim}{=}}

\newcommand{\shHom}{\underline{\operatorname{Hom}}}

\newcommand{\shAut}{\underline{\operatorname{Aut}}}
\newcommand{\Ext}{\operatorname{Ext}}

\newcommand{\Leib}{\operatorname{Leib}}
\newcommand{\Lie}{\operatorname{Lie}}

\newcommand{\CA}{{\mathcal CA}}
\newcommand{\ECA}{{\mathcal ECA}}
\newcommand{\VA}{{\mathcal VA}}
\newcommand{\EVA}{{\mathcal EVA}}
\newcommand{\DR}{{\Omega^\bullet_X}}

\newcommand{\vac}{{\mathbf{1}}}
\begin{document}

\title{Vertex Algebroids I}

\author{Paul Bressler}
\address{}
\email{bressler@@math.arizona.edu}
\date{\today}
\maketitle

\section{Introduction}
The purpose of this note is to give a ``coordinate free'' construction
and prove the uniqueness of the vertex algebroid which gives rise to
the chiral de Rham complex of \cite{GMS}. In order to do that, we adapt
the strategy of \cite{BD} to the setting of vertex agebroids.

To this end we show that the stack $\EVA_X$ of exact vertex algebroids on $X$
is a torsor under the stack in Picard groupoids $\ECA_X$ of exact Courant algebroids
on $X$ (and, in particular, is locally non-empty). Moreover, we show that that
$\ECA_X$ is naturally equivalent to the stack of
$\Omega^2_X @>d>> \Omega^{3,cl}_X$-torsors. These facts are proven for $X$ a manifold
($C^\infty$, (complex) analytic, smooth algebraic variety over $\bbC$). We leave
to the reader the obvious extension to differential graded manifolds.

Given a manifold $X$, let $X^\sharp$ denote the differential graded manifold with the
underlying space $X$ and the structure sheaf $\cO_{X^\sharp}$ the de Rham complex of
$X$. We show that $\ECA_{X^\sharp}$ is, in fact, trivial. Together with local existence
of vertex $\cO_{X^\sharp}$-algebroids this implies that there exists a unique up to a
unique isomorphism vertex algebroid over the de Rham complex. We give a ``coordinate-free''
description of the unique vertex $\cO_{X^\sharp}$-algebroid in terms of
generators and relations.

\section{Vertex operator algebras}
Throughout this section we follow the notations of \cite{GMS}. The following
definitions are lifted from loc. cit.

\subsection{}
A {\it $\bbZ_{\geq 0}$-graded vertex algebra} is a
$\bbZ_{\geq 0}$-graded $k$-module $V=\oplus\ V_i$, equipped with
a distinguished vector $\vac\in V_0$ ({\it vacuum vector})
and a family of bilinear operations
\[
_{(n)}:\ V\times V @>>> V,\ (a, b)\mapsto a_{(n)}b
\]
of degree $-n-1$, $n\in \bbZ$, such that
\[
\vac_{(n)}a=\delta_{n,-1}a;\ a_{(-1)}\vac=a;\ a_{(n)}\vac=0\
\text{\ if\ }n\geq 0,
\]
and
\[
\sum_{j=0}^\infty\ \binom{m}{j}(a_{(n+j)}b)_{(m+l-j)}c= \\
\sum_{j=0}^\infty\ (-1)^j\binom{n}{j}
\bigl\{a_{(m+n-j)}b_{(l+j)}c-(-1)^nb_{(n+l-j)}a_{(m+j)}c\bigr\}
\]
for all $a,b,c\in V,\ m,n,l\in \bbZ$.

\subsection{}
A morphism of vertex algebras is a map of graded $k$-modules (of degree zero)
which maps the vacuum vector to the vacuum vector and commutes with all of
the operations.

Let $\cV ert$ denote the category of vertex algebras.

\subsection{}
Let
\[
\partial^{(j)}a:=a_{(-1-j)}\b1,\ \ j\in\bbZ_{\geq 0} \ .
\]
Then, $\partial^{(j)}$ is an endomorphisms of $V$ of degree $j$ which satisfies
(see \cite{GMS})
\begin{itemize}
\item $\partial^{(j)}\vac=\delta_{j,0}\vac$,
\item $\partial^{(0)}=Id$,
\item $\partial^{(i)}\cdot\partial^{(j)}=\binom{i+j}{i}\partial^{(i+j)}$,
\item $(\partial^{(j)}a)_{(n)}b=(-1)^j\binom{n}{j}a_{(n-j)}b$, and
\item $\partial^{(j)}(a_{(n)}b)=\sum_{p=0}^j\ (\partial^{(p)}a)_{(n)}\partial^{(j-p)}b$
\end{itemize}
for all $n\in\bbZ$.

\subsection{}
The subject of the definition below is the restriction of the structure of a vertex
algebra to the graded components of degrees zero and one.

A {\it $1$-truncated vertex algebra} is a
septuple $v=(V_0,V_1,\vac,\partial,_{(-1)},_{(0)},_{(1)})$ where
\begin{itemize}
\item $V_0, V_1$ are $k$-modules,
\item $\vac$ an element of $V_0$ ({\it vacuum vector}),
\item $\partial:\ V_0 @>>> V_1$ a $k$-linear map,
\item $_{(i)}:\ (V_0\oplus V_1)\times (V_0\oplus V_1) @>>> V_0\oplus V_1$
(where $(i=-1,0,1)$) are $k$-bilinear operations of degree $-i-1$.
\end{itemize}

Elements of $V_0$
(resp., $V_1$) will be denoted $a,b,c$ (resp., $x,y,z$). There are
seven operations: $a_{(-1)}b, a_{(-1)}x, x_{(-1)}a, a_{(0)}x$, $
x_{(0)}a, x_{(0)}y$ and $x_{(1)}y$. These satisfy the following axioms:
\begin{itemize}
\item
(Vacuum)
\[
a_{(-1)}\vac=a;\ x_{(-1)}\vac=x;\ x_{(0)}\vac=0
\]
\item
(Derivation)
\begin{itemize}
\item[$Deriv_1\ \ $]
$(\partial a)_{(0)}b=0;\ (\partial a)_{(0)}x=0;\
(\partial a)_{(1)}x=-a_{(0)}x$
\item[$Deriv_2\ \ $]
$\partial(a_{(-1)}b)=(\partial a)_{(-1)}b+a_{(-1)}\partial b;\
\partial(x_{(0)}a)=x_{(0)}\partial a$
\end{itemize}
\item
(Commutativity)
\begin{itemize}
\item[$Comm_{-1}$]
$a_{(-1)}b=b_{(-1)}a;\ a_{(-1)}x=x_{(-1)}a-\partial(x_{(0)}a)$
\item[$Comm_0\ \ $]
$x_{(0)}a=-a_{(0)}x;\ x_{(0)}y=-y_{(0)}x+\partial(y_{(1)}x)$
\item[$Comm_1\ \ $]
$x_{(1)}y=y_{(1)}x$
\end{itemize}
\item
(Associativity)
\begin{itemize}
\item[$Assoc_{-1}$] $(a_{(-1)}b)_{(-1)}c=a_{(-1)}b_{(-1)}c$
\item[$Assoc_0\ \ $]
$\alpha_{(0)}\beta_{(i)}\gamma=
(\alpha_{(0)}\beta)_{(i)}\gamma+\beta_{(i)}\alpha_{(0)}\gamma,\ (\alpha,
\beta, \gamma\in V_0\oplus V_1)$
whenever the both sides are defined, i.e. the operation $_{(0)}$ is a derivation
of all of the operations $_{(i)}$.
\item[$Assoc_1\ \ $]
$(a_{(-1)}x)_{(0)}b=a_{(-1)}x_{(0)}b$
\item[$Assoc_2\ \ $]
$(a_{(-1)}b)_{(-1)}x=a_{(-1)}b_{(-1)}x+(\partial a)_{(-1)}b_{(0)}x+
(\partial b)_{(-1)}a_{(0)}x$
\item[$Assoc_3\ \ $]
$(a_{(-1)}x)_{(1)}y=a_{(-1)}x_{(1)}y-x_{(0)}y_{(0)}a$
\end{itemize}
\end{itemize}

\subsection{}
A {\it morphism} between two $1$-truncated vertex algebras
$f:\ v=(V_0,V_1,\ldots) @>>> v'=(V'_0,V'_1,\ldots)$ is a pair of maps of $k$-modules
$f=(f_0,f_1),\ f_i:\ V_i @>>> V'_i$ such that $f_0(\vac)=\vac',\ f_1(\partial a)=
\partial f_0(a)$ and $f(\alpha_{(i)}\beta)=f(\alpha)_{(i)}f(\beta)$,
whenever both sides are defined.

Let $\cV ert_{\leq 1}$ denote the category of $1$-truncated vertex algebras.
We have an obvious truncation functor
$$
t:\ \cV ert @>>> \cV ert_{\leq 1}
\eqno{(3.1.1)}
$$
which assignes to a vertex algebra $V=\oplus V_i$ the truncated algebra
$tV:=(V_0,V_1,\ldots)$.

\subsection{Remark}\label{rem:comm-alg}
It follows easily that the operation $_{(-1)}\ :V_0\times V_0 @>>> V_0$
endows $V_0$ with a structure of a commutative $k$-algebra.

\section{Vertex algebroids}
Suppose that $X$ is smooth variety over $\bbC$ (a complex manifold,
a $C^\infty$-manifold). In either case we will denoted by $\cO_X$
($\cT_X$, $\Omega^i_X$) the corresponding structure sheaf (the sheaf
of vector fields, the sheaf of differential $i$-forms).

A vertex $\cO_X$-algebroid, as defined in this section, is, essentially,
a sheaf of $1$-truncated vertex algebras, whose degree zero component (which is
a sheaf of algebras by \ref{rem:comm-alg}) is identified with $\cO_X$.

\subsection{Leibniz algebras}
A {\em Leibnitz $k$-algebra} is a $k$-vector space $\frak g$ equipped with
a bilinear operation $[\ ,\ ]:\frak g\otimes g @>>> g$ (the Leibniz bracket)
which satisfies the Jacobi type identity $[a,[b,c]] = [[a,b],c] +
[b,[a,c]]$.

A morphism of Leibniz $k$-algebras is a $k$-linear map which commutes with
the respective Leibniz brackets.

We denote the category of Leibniz $k$-algebras by $\Leib(k)$.
Let $\Lie(k)$ denote the category of Lie $k$-algebras.

Since any Lie algebra is a Leibnitz algebra there is an full
embedding
\begin{equation}\label{lie-leib}
\Lie(k) @>>> \Leib(k) \ .
\end{equation}
%A Leibnitz algebra $\frak g$ has the maximal Lie quotient
%${\frak g}_{Lie}$ by the relations $[a,b]+[b,a]=0$, $a,b\in\frak g$.
%The assignment ${\frak g}\mapsto{\frak g}_{Lie}$ extends to the
%functor
%\[
%(\ )_{Lie} : \Leib(k) @>>> \Lie(k)
%\]
%which is right adjoint to \eqref{lie-leib}.

\subsection{Vertex algebroids}
A {\em vertex $\cO_X$-algebroid} is a sheaf of $\bbC$-vector
spaces $\cV$ with a pairing
\begin{eqnarray*}
\cO_X\otimes_\bbC\cV & @>>> & \cV \\
f\otimes v & \mapsto & f*v
\end{eqnarray*}
such that $1* v = v$ (i.e. a ``non-associative unital
$\cO_X$-module'') equipped with
\begin{enumerate}
\item
a structure of a Leibniz $\bbC$-algebra $[\ ,\ ] :
\cV\otimes_\bbC\cV @>>> \cV$

\item
a $\bbC$-linear map of Leibniz algebras $\pi : \cV @>>> \cT_X$
(the {\em anchor})
\item
a symmetric $\bbC$-bilinear pairing $\langle\ ,\ \rangle :
\cV\otimes_\bbC\cV @>>> \cO_X$
\item
a $\bbC$-linear map $\partial : \cO_X @>>> \cV$ such that
$\pi\circ\partial = 0$
\end{enumerate}
which satisfy
\begin{eqnarray}
f*(g*v) - (fg)*v & = & - \pi(v)(f)*\partial(g) -
\pi(v)(g)*\partial(f)\label{assoc} \\
\left[v_1,f*v_2\right] & = & \pi(v_1)(f)*v_2 + f*[v_1,v_2] \label{leib}
\\
\left[v_1,v_2\right] + [v_2,v_1] & = & \partial(\langle v_1,v_2\rangle)
\label{symm-bracket}\\
\pi(f*v) & = & f\pi(v) \label{anchor-lin} \\
\langle f*v_1, v_2\rangle & = & f\langle v_1,v_2\rangle -
\pi(v_1)(\pi(v_2)(f)) \label{pairing}\\
\pi(v)(\langle v_1, v_2\rangle) & = & \langle[v,v_1],v_2\rangle +
\langle v_1,[v,v_2]\rangle \label{pairing-inv} \\
\partial(fg) & = & f*\partial(g) + g*\partial(f) \label{deriv} \\
\left[v,\partial(f)\right] & = & \partial(\pi(v)(f)) \label{bracket-o}\\
\langle v,\partial(f)\rangle & = & \pi(v)(f)\label{pairing-o}
\end{eqnarray}
for $v,v_1,v_2\in\cV$, $f,g\in\cO_X$.

\subsubsection{}
A morphism of vertex $\cO_X$-algebroids is a $\bbC$-linear map of
sheaves which preserves all of the structures.

We denote the category of vertex $\cO_X$-algebroids by
$\VA_{\cO_X}(X)$. It is clear that the notion vertex $\cO_X$-algebroid
is local, i.e. vertex $\cO_X$-algebroids form a stack which we denote
by $\VA_{\cO_X}$.

\subsection{From vertex algebroids to 1-truncated vertex algebras}
Suppose that $\cV$ is a vertex $\cO_X$-algebroid. For $f,g\in\cO_X$,
$v,w\in\cV$ let
\begin{equation}\label{op:-1}
f_{(-1)}g = fg,\ \ f_{(-1)}v = f*v,\ \ v_{(-1)}f = f*v - \partial\pi(v)(f),
\end{equation}
\begin{equation}\label{op:0}
v_{(0)}f = - f_{(0)}v = \pi(v)(f),\ \ v_{(0)}w = [v,w],
\end{equation}
\begin{equation}\label{op:1}
v_{(1)}w = \langle v,w\rangle\ .
\end{equation}

\subsubsection{Lemma}
The septuple $(\cO_X,\cV,1,\partial, _{(-1)}, _{(0)}, _{(1)})$ is a
sheaf of 1-truncated vertex operator algebras.

Conversely, if the septuple $(\cO_X,\cV,1,\partial,_{(-1)},_{(0)},_{(1)})$
is a sheaf of 1-truncated vertex operator algebras, then the formulas
\eqref{op:-1}, \eqref{op:0}, \eqref{op:1} define a structure of a vertex
$\cO_X$-algebroid on $\cV$.

\subsection{The associated Lie algebroid}
Suppose that $\cV$ is a vertex $\cO_X$-algebroid. Let
\begin{eqnarray*}
\Omega_\cV & \stackrel{def}{=} & \cO_X*\partial(\cO_X)\subset\cV \ ,\\
\overline\cV & \stackrel{def}{=} & \cV/\Omega_\cV \ .
\end{eqnarray*}
Note that the symmetrization of the Leibniz bracket takes values
in $\Omega_\cV$.

For $f,g,h\in\cO_X$
\[
f*(g*\partial(h)) - (fg)*\partial(h) =
\pi(\partial(h))(f)*\partial(g) + \pi(\partial(h))(g)*\partial(f)
= 0\ ,
\]
because $\pi\circ\partial = 0$. Therefore, $\cO_X*\Omega_\cV =
\Omega_\cV$, and $\Omega_\cV$ is an $\cO_X$-module. The map
$\partial : \cO_X @>>> \Omega_\cV$ is a derivation, hence induces
the $\cO_X$-linear map $\Omega^1_X @>>> \Omega_\cV$.

Since the associator of the $\cO_X$-action on $\cV$ takes values
in $\Omega_\cV$, $\overline\cV$ is an $\cO_X$-module.

For $f,g,h\in\cO_X$
\[
\pi(f\partial(g))(h) = f\pi(\partial(g))(h) = 0 \ .
\]
Therefore, $\pi$ vanishes on $\Omega_\cV$, hence, factors through the
map
\begin{equation}\label{VA-anchor}
\pi : \overline\cV @>>> \cT_X
\end{equation}
of $\cO_X$-modules.

For $v\in\cV$, $f,g\in\cO_X$
\begin{eqnarray*}
[v,f\partial(g)] & = & \pi(v)(f)\partial(g)+f[v,\partial(g)] \\
& = & \pi(v)(f)\partial(g)+f\partial(\pi(v)(g)) \ .
\end{eqnarray*}
Therefore, $[\cV,\Omega_\cV]\subseteq\Omega_\cV$ and the Leibniz bracket on
$\cV$ descends to the operation
\begin{equation}\label{VA-bracket}
[\ ,\ ]:\overline\cV\otimes_\bbC\overline\cV @>>> \overline\cV
\end{equation}
which is skew-symmetric because the symmetrization of the Leibniz
bracket on $\cV$ takes values in $\Omega_\cV$ and satisfies the Jacobi
identity because the Leibniz bracket on $\cV$ does.

\subsubsection{Lemma}
The $\cO_X$-module $\overline\cV$ with the bracket
\eqref{VA-bracket} and the anchor \eqref{VA-anchor} is a Lie
$\cO_X$-algebroid.

\subsection{Transitive vertex algebroids}
A vertex $\cO_X$-algebroid is called {\em transitive} if the anchor
map is surjective.

\subsubsection{Remark}
The vertex $\cO_X$-algebroid $\cV$ is called transitive if and only if
the Lie $\cO_X$-algebroid $\overline\cV$ is.

\subsubsection{}
Suppose that $\cV$ is a transitive vertex $\cO_X$-algebroid. The
derivation $\partial$ induces the map
\[
i : \Omega^1_X @>>> \cV \ .
\]
For $v\in\cV$, $f,g\in\cO_X$
\begin{eqnarray*}
\langle v,f\partial(g)\rangle & = & f\langle v,\partial(g)\rangle
- \pi(\partial(g))\pi(v)(f) \\
& = & f\pi(v)(g) \\
& = & \iota_{\pi(v)}fdg \ .
\end{eqnarray*}
If follows that the map $i$ is adjoint to the anchor map $\pi$.
The surjectivity of the latter implies the latter implies the
injectivity of the former. Since, in addition, $\pi\circ i = 0$
the sequence
\[
0 @>>> \Omega^1_X @>i>> \cV @>\pi>> \overline\cV @>>> 0
\]
is exact and $i$ is isotropic.

\subsection{Exact vertex algebroids}
A vertex algebroid $\cV$ is called {\em exact} if the map
$\overline\cV @>>> \cT_X$ is an isomorphism.

We denote the stack of exact vertex $\cO_X$-algebroids by
$\EVA_X$.

A morphism of exact vertex algebroids induces a morphism of
respective extensions of $\cT_X$ by $\Omega^1_X$, hence is an
isomorphism of sheaves of $\bbC$-vector spaces. It is clear that
the inverse isomorphism is a morphism of vertex
$\cO_X$-algebroids. Hence, $\EVA_X$ is a stack in groupoids.

\subsubsection{Example}\label{example:loc-pic}
Suppose that $\cT_X$ is freely generated as an $\cO_X$-module by
a locally constant subsheaf of Lie $\bbC$-subalgebras
$\tau\subset\cT_X$, i.e. the canonical map $\cO_X\otimes_\bbC\tau
@>>> \cT_X$ is an isomorphism.

There is a unique structure of an exact vertex $\cO_X$-algebroid on
$\cV=\Omega^1_X\bigoplus\cT_X$ such that
\begin{itemize}
\item $f*(1\otimes t) = f\otimes t$ for $f\in\cO_X$, $t\in\tau$,

\item the anchor map is given by the projection $\cV @>>> \cT_X$,

\item the map $\tau @>>> \cV$ is a morphism of Leibniz algebras,

\item the pairing on $\cV$ restricts to the trivial pairing on
$\tau$,

\item the derivation $\partial : \cO_X @>>> \cV$ is given by the
composition $\cO_X @>d>>\Omega^1_X @>>> \cV$.
\end{itemize}

Indeed, the action of $\cO_X$ is completely determined by \eqref{assoc}:
for $f,g\in\cO_X$, $t\in\tau$,
\[
f*(g\otimes t) = f*(g*(1\otimes t)) = fg\otimes t - t(f)dg - t(g)df \ .
\]
In a similar fashion the bracket is completely determined by \eqref{leib}
and \eqref{symm-bracket}, and the pairing is determined by \eqref{pairing}.

We leave the verification of the identities \eqref{assoc} - \eqref{pairing-o}
to the reader.

\section{Courant algebroids}
Courant algebroids are classical limits of vertex algebroids. They are
related to vertex Poisson algebras (coisson algebras in the terminology of
\cite{BD}) in the same way as the vertex algebroids
are related to vertex operators algebras.

\subsection{Courant algebroids}
A {\em Courant $\cO_X$-algebroid} is an $\cO_X$-module $\cQ$
equipped with
\begin{enumerate}
\item a structure of a Leibniz $\bbC$-algebra
\[
[\ ,\ ] : \cQ\otimes_\bbC\cQ @>>> \cQ \ ,
\]

\item
an $\cO_X$-linear map of Leibniz algebras (the anchor map)
\[
\pi : \cQ @>>> \cT_X \ ,
\]

\item
a symmetric $\cO_X$-bilinear pairing
\[
\langle\ ,\ \rangle : \cQ\otimes_{\cO_X}\cQ @>>> \cO_X \ ,
\]

\item
a derivation
\[
\partial : \cO_X @>>> \cQ
\]
such that $\pi\circ\partial = 0$
\end{enumerate}
which satisfy
\begin{eqnarray}
[q_1,fq_2] & = & f[q_1,q_2] + \pi(q_1)(f)q_2 \\
\langle [q,q_1],q_2\rangle + \langle q_1,[q,q_2]\rangle & = & \pi(q)(\langle q_1, q_2\rangle) \\
\left[q,\partial(f)\right] & = & \partial(\pi(q)(f)) \\
\langle q,\partial(f)\rangle & = & \pi(q)(f)\label{axiom:courant-ip-with-df} \\
\left[q_1,q_2\right] + [q_2,q_1] & = & \partial(\langle q_1, q_2\rangle)
\end{eqnarray}
for $f\in\cO_X$ and $q,q_1,q_2\in\cQ$.

\subsubsection{}
A morphism of Courant $\cO_X$-algebroids is an $\cO_X$-linear map
of Leibnitz algebras which commutes with the respective anchor
maps and derivations and preserves the respective pairings.

We denote the category of Courant $\cO_X$-algebroids on $X$ by
$\CA_{\cO_X}(X)$. The notion of Courant $\cO_X$-algebroid is local, i.e.
Courant $\cO_X$-algebroids form a stack which we denote $\CA_{\cO_X}$.

\subsection{The associated Lie algebroid}
Suppose that $\cQ$ is a Courant $\cO_X$-algebroid. Let
\begin{eqnarray*}
\Omega_\cQ & \stackrel{def}{=} & \cO_X\partial(\cO_X)\subset\cQ \ , \\
\overline\cQ & \stackrel{def}{=} & \cQ/\Omega_\cQ \ .
\end{eqnarray*}
Note that the symmetrization of the Leibniz bracket on $\cQ$ takes
values in $\Omega_\cQ$.

For $q\in\cQ$, $f,g\in\cO_X$
\begin{eqnarray*}
[q,f\partial(g)] & = & f[q,\partial(g)] + \pi(q)(f)\partial(g) \\
& = & f\partial(\pi(q)(g)) + \pi(q)(f)\partial(g)
\end{eqnarray*}
which shows that $[\cQ,\Omega_\cQ]\subseteq\Omega_\cQ$. Therefore the
Leibniz bracket on $\cQ$ descends to the Lie bracket
\begin{equation}\label{LAbracket}
[\ ,\ ] : \overline\cQ\otimes_\bbC\overline\cQ @>>> \overline\cQ\ .
\end{equation}

Because $\pi$ is $\cO_X$-linear and $\pi\circ\partial = 0$, $\pi$
vanishes on $\Omega_\cQ$ and factors through the map
\begin{equation}\label{LAanchor}
\pi : \overline\cQ @>>> \cT_X \ .
\end{equation}

\subsubsection{Lemma}
The bracket \eqref{LAbracket} and the anchor \eqref{LAanchor}
determine the structure of a Lie $\cO_X$-algebroid on
$\overline\cQ$.

\subsection{Transitive Courant algebroids}
A Courant $\cO_X$-algebroid is called {\em transitive} if the
anchor map is surjective.

\subsubsection{Remark}
A Courant $\cO_X$-algebroid $\cQ$ is transitive if and only if the
associated Lie $\cO_X$-algebroid is.

\subsubsection{}
Suppose  that  $\cQ$ is a transitive Courant $\cO_X$-algebroid. The
derivation $\partial$ induces the $\cO_X$-linear map
\[
i : \Omega^1_X @>>> \cQ \ .
\]
Since $\langle q, \alpha\rangle = \iota_{\pi(q)}\alpha$, it
follows that the map $i$ is adjoint to the anchor map $\pi$. The
surjectivity of the latter implies that $i$ is injective. Since,
in addition, $\pi\circ i = 0$ the sequence
\[
0 @>>> \Omega^1_X @>{i}>> \cQ @>>> \overline\cQ @>>> 0
\]
is exact. Moreover, $i$ is isotropic with respect to the
symmetric pairing.

\subsection{Exact Courant algebroids}
The Courant algebroid $\cQ$ is called {\em exact} if the anchor
map $\pi: \overline{\cQ} @>>> \cT_X$ is an isomorhism.

We denote the stack of exact Courant $\cO_X$-algebroids by
$\ECA_X$.

A morphism of exact Courant algebroids induces a morphism of
respective extensions of $\cT_X$ by $\Omega^1_X$, hence is an
isomorphism of $\cO_X$-modules. It is clear that the inverse is a
morphism of Courant $\cO_X$-algebroids. Hence $\ECA_X$ is
a stack in groupoids.

\subsection{Automorphisms of exact Courant algebroids}
Let $\Ext^{\langle\ ,\ \rangle}_{\cO_X}(\cT_X,\Omega^1_X)$ denote
the category whose objects are pairs $(\cE,\langle\ ,\ \rangle)$,
where $\cE$ is an extension
\[
0 @>>> \Omega^1_X @>i>> \cE @>\pi>> \cT_X @>>> 0 \ .
\]
and
\[
\langle\ ,\ \rangle : \cE\otimes_{\cO_X}\cE @>>> \cO_X
\]
is a symmetric pairing such that $i$ is Lagrangian and the induced
pairing between $\Omega^1_X$ and $\cT_X=\cE/\Omega^1_X$ is the
canonical duality pairing. A morphism of such is a morphism of the
underlying extensions which preserves the respective pairings.

The map
\begin{equation}\label{Ext-exp}
\exp : \shHom_{\cO_X}(\cT_X,\Omega^1_X) @>>> \shAut_{\Ext}(\cE)
\end{equation}
defined by $\phi\mapsto (e\mapsto e + \phi(\pi(e)))$ is an
isomorphism. It restricts to the isomorphism
\[
\exp : \Omega^2_X @>>> \shAut_{\Ext}^{\langle\ ,\
\rangle}(\cE,\langle\ ,\ \rangle) \ ,
\]
via the map
\[
\Omega^2_X @>>> \shHom_{\cO_X}(\cT_X,\Omega^1_X)
\]
defined by $\beta\mapsto(\xi\mapsto\iota_\xi\beta)$.

Suppose that $\cQ$ is an exact Courant $\cO_X$-algebroid. The
automorphism, induced by a 2-form $\beta$ of the underlying
extension preserves the Leibnitz bracket if and only if the form
$\beta$ is closed, i.e. the map \eqref{Ext-exp} restricts to the
isomorphism
\[
\Omega^{2,cl}_X @>>> \shAut_\ECA(\cQ) \ .
\]

\subsection{The $\bbC$-vector space structure of $\ECA$}
The category $\ECA_X(X)$ has a natural structure of a
``$\bbC$-vector space in categories'' induced by that of
$\Ext^1_{\cO_X}(\cT_X,\Omega^1_X)$.

\subsubsection{Addition}
Suppose that $\cQ_1$ and $\cQ_2$ are two exact Courant
$\cO_X$-algebroids. Let $\cQ_1+\cQ_2$ denote the push-out of
$\cQ_1\times_{\cT_X}\cQ_2$ by the addition map
$\Omega^1_X\times\Omega^1_X @>>> \Omega^1_X$. A section of
$\cQ_1+\cQ_2$ is an equivalence class of pairs $(q_1,q_2)$, where
$q_i\in\cQ_i$ and $\pi(q_1)=\pi(q_2)$. Two pairs are equivalent if
their (componentwise) difference is of the form $(i(\alpha),
-i(\alpha))$ for some $\alpha\in\Omega^1_X$.

The two maps $\Omega^1_X @>>> \cQ_1+\cQ_2$ given by
$\alpha\mapsto(i(\alpha),0)$ and $\alpha\mapsto(0,i(\alpha))$
coincide. We denote their common value by $i$ as well. There is a
short exact sequence
\[
0 @>>> \Omega^1_X @>i>> \cQ_1+\cQ_2 @>\pi>> \cT_X @>>> 0 \ ,
\]
where $\pi((q_1,q_2))$ is defined as the common value of
$\pi(q_1)$ and $\pi(q_2)$.

The map $i$ determines the derivation $\partial :\cO_X @>>>
\cQ_1+\cQ_2$ by $\partial(f) = i(df)$.

For $q_i,q'_i\in\cQ_i$ let
\begin{equation}\label{formula:sum-bracket-pairing}
[(q_1,q_2),(q'_1,q'_2)] = ([q_1,q'_1],[q_2,q'_2]), \ \ \ \langle
(q_1,q_2),(q'_1,q'_2)\rangle = \langle q_1,q'_1\rangle + \langle
q_2,q'_2\rangle \ .
\end{equation}

\subsubsection{Lemma}
The bracket and the symmetric pairing given by
\eqref{formula:sum-bracket-pairing} determine a structure of an
exact Courant algebroid on $\cQ_1+\cQ_2$.

\subsubsection{Scalar multiplication}
Suppose that $\cQ$ is an exact Courant $\cO_X$-algebroid and
$\lambda\in\bbC$. Let $\lambda\cQ$ denote push-out of $\cQ$ by the
multiplication by $\lambda$ map $\Omega^1_X
@>{\lambda\cdot}>>\Omega^1_X$. A section of $\lambda\cQ$ is an
equivalence class of pairs $(\alpha,q)$ with $\alpha\in\Omega^1_X$
and $q\in\cQ$. Two pairs as above are equivalent if their
componentwise difference is of the form
$(\lambda\alpha,-i(\alpha))$ for some $\alpha\in\Omega^1_X$.

Let $i:\Omega^1_X @>>> \lambda\cQ$ denote the map
$\alpha\mapsto(\alpha,0)$. There is a short exact sequence
\[
0 @>>> \Omega^1_X @>i>> \lambda\cQ @>\pi>> \cT_X @>>> 0
\]
where $\pi(\alpha,q) = \pi(q)$.

The map $i$ determines the derivation $\partial:\cO_X @>>>
\lambda\cQ$ by $\partial(f) = i(df) = (df,0)$. Note that $\pi\circ
i =0$ holds.

\subsubsection{Lemma}
There is a unique structure of exact Courant algebroid on
$\lambda\cQ$ (with anchor $\pi$ and derivation $\partial$ as
above) such that the map $\cQ @>>> \lambda\cQ: q\mapsto(0,q)$ is a
morphism of Leibniz algebras.

\begin{pf}
For $q_i\in\cQ$ the calculation
\begin{multline*}
(\lambda d\langle q_1,q_2\rangle,0)= (0,\partial\langle
q_1,q_2\rangle) = \\
(0,[q_1,q_2] + [q_2,q_1]) = [(0,q_1),(0,q_2)]
+ [(0,q_2),(0,q_1)] =
\partial\langle(0,q_1),(0,q_2)\rangle = \\
(d\langle(0,q_1),(0,q_2)\rangle,0)
\end{multline*}
together with \eqref{axiom:courant-ip-with-df} determines the
symmetric pairing on $\lambda\cQ$. In particular, $
\langle(0,q_1),(0,q_2)\rangle=\lambda\langle q_1,q_2\rangle$.
\end{pf}

\subsubsection{}
Let $\cQ_0$ denote the exact Courant algebroid with the underlying
$\cO_X$-module the trivial extension $\Omega^1_X\oplus\cT_X$,
self-evident derivation and anchor map, the symmetric paring the
obvious extension of the natural duality pairing and the Leibniz
bracket uniquely determined by the requirement that the natural
inclusion of $\cT_X$ is a map of Leibniz algebras.

Note that $\cQ_0$ is equipped with a flat connection. For any
exact Courant algebroid $\cQ$ a flat connection on $\cQ$ is the
same thing as a morphism $\cQ_0 @>>> \cQ$.

\subsubsection{Lemma}
For any exact Courant algebroid $\cQ$ there is a canonical
isomorphism $\cQ\isomo\cQ+\cQ_0$.

\subsubsection{Remark}
The natural morphism (of Leibniz algebras) $\cQ @>>> \lambda\cQ$
factors through $\cT_X$ if $\lambda=0$ (and is an isomorphism of
the underlying $\cO_X$-modules otherwise) and the composition
$\cT_X @>>> 0\cQ @>>> \cT_X$ is clearly equal to the identity.
Thus, $0\cQ$ carries a canonical flat connection, or,
equivalently, is canonically isomorphic to $\cQ_0$.

\subsubsection{Remark}
If $\lambda\neq 0$ the natural map $\cQ @>>> \lambda\cQ$ is an isomorphism
of the underlying $\cO_X$-modules with the inverse given by
$(\alpha,q)\mapsto q+\lambda^{-1}\alpha$. The Courant algebroid structure
on $\cQ$ induced via this isomorphism and denoted $\partial_{\lambda\cQ}$,
etc. is given by
\[
\partial_{\lambda\cQ}(f) = \lambda^{-1}\partial(f),
 \ \ \ \langle q_1,q_2\rangle_{\lambda\cQ} = \lambda\langle q_1,q_2\rangle,
 \ \ \ [q_1,q_2]_{\lambda\cQ}=[q_1,q_2]
\]

\subsection{Classification of exact Courant algebroids}
\subsubsection{Connections}
A {\em connection} on an exact Courant $\cO_X$-algebroid $\cQ$ is
a Lagrangian section of the anchor map.

\subsubsection{Curvatrue}\label{ECAcurv}
Suppose that $\cQ$ is an exact Courant $\cO_X$-algebroid and
$\nabla$ is a connection on $\cQ$. The formula
\[
(\xi,\xi_1,\xi_2)\mapsto\iota_\xi([\nabla(\xi_1),\nabla(\xi_2)] -
\nabla([\xi_1,\xi_2])
\]
defines a differential 3-form $c(\nabla)$ on $X$ which is easily
seen to be closed.

The differential form $c(\nabla)$ is
called {\em the curvature} of the connection $\nabla$.

\subsubsection{}
Suppose that $\cQ$ is an exact Courant $\cO_X$-algebroid, $\nabla$
is a connection on $\cQ$ and
$\phi\in\Hom_{\cO_X}(\cT_X,\Omega^1_X)$. Then, $\nabla+\phi$ is
another section of the anchor map. The section $\nabla+\phi$ is
Lagrangian if and only if $\phi$ is derived from a differential
2-form, i.e. there exists $\beta\in\Omega^2_X(X)$ such that
$\phi(\xi) = \iota_\xi\beta$ for all $\xi$.

Let $\cC(\cQ)$ denote the sheaf of (locally defined) connections
on $\cQ$. It is a $\Omega^2_X$-torsor. The correspondence
$\nabla\mapsto c(\nabla)$ defines the map
\[
c : \cC(\cQ) @>>> \Omega^{3,cl}_X
\]
which satisfies
\[
c(\nabla + \widetilde\beta) = c(\nabla) + d\beta \ ,
\]
where $\widetilde{(\ )} : \Omega^i_X @>>>
\shHom_{\cO_X}(\cT_X,\Omega^{i-1}_X)$ is defined by
$\widetilde\beta(\xi) = \iota_\xi\beta$.

Therefore, the pair $(\cC(\ ),c)$ defines the functor
\begin{equation}\label{ECA-to-tors}
\ECA_X(X) @>>> (\Omega^2_X @>d>> \Omega^{3,cl})-tors \ .
\end{equation}

\subsubsection{Lemma}
The functor \eqref{ECA-to-tors} is an eqivalence of $\bbC$-vector
spaces in categories.

\subsubsection{Corollary}
The $\bbC$-vector space of isomorphism classes of exact Courant
$\cO_X$-algebroids is naturally isomorphic to
$H^1(X;\Omega^2_X @>d>> \Omega^{3,cl}_X)$.

\subsection{The action on $\EVA$}
The groupoid $\EVA_X(X)$, if non-empty, is an ``affine space in
categories'' (under the action of the ``$\bbC$-vector space in categories''
$\ECA_X(X)$).

\subsubsection{}
Suppose that $\cV$ (respectively, $\cQ$) is an exact vertex (respectively,
Courant) algebroid. Let $\cQ+\cV$ denote the push-out of
$\cQ\times_{\cT_X}\cV$ by the addition map $\Omega^1_X\times\Omega^1_X @>+>>
\Omega^1_X$. A section of $\cQ+\cV$ is an equivalence class of pairs $(q,v)$
with $q\in\cQ$, $v\in\cV$ and $\pi(q)=\pi(v)$. Two pairs are equivalent if
their (componentwise) difference is of the form $(i(\alpha),-i(\alpha))$
for some $\alpha\in\Omega^1_X$.

The two maps $\Omega^1_X @>>> \cQ+\cV$ given by $\alpha\mapsto(i(\alpha),0)$
and $\alpha\mapsto(0,i(\alpha))$ coincide and we denote their common value
by $i$ as well. There is a short exact sequence
\[
0 @>>> \Omega^1_X @>i>> \cQ+\cV @>\pi>> \cT_X @>>> 0 \ ,
\]
where $\pi((q,v))$ is defined to be the common value of $\pi(q)$ and $\pi(v)$.
Let
\begin{multline}\label{formulas:Q-plus-V}
\partial(f)=i(df),\ \ \ f*(q,v)=(fq,f*v), \\ 
[(q_1,v_1),(q_2,v_2)]= ([q_1,q_2],[v_1,v_2]),
\ \ \ \langle(q_1,v_1),(q_2,v_2)\rangle =
\langle q_1,q_2\rangle+\langle v_1,v_2\rangle
\end{multline}

\subsubsection{Lemma}
The formulas \eqref{formulas:Q-plus-V} determine a structure of
a vertex algebroid on $\cQ+\cV$.

\subsubsection{Proposition}
$\EVA_X$ is a torsor under $\ECA_X$.
\begin{pf}
Suppose that $\cV_1$ and $\cV_2$ are two exact vertex $\cO_X$-algebroids.
Let $\cV_2-\cV_1$ the push-out of $\cV_2\times_{\cT_X}\cV_1$ by the difference
map $\Omega^1_X\times\Omega^1_X @>->>\Omega^1_X:(\alpha,\beta)\mapsto\alpha-\beta$.
A section of $\cV_2-\cV_1$ is an equivalence class of pairs $(v_1,v_2)$, where
$v_i\in\cV_i$ and $\pi(v_1)=\pi(v_2)$. Two pairs are equivalent if
their (componentwise) difference is of the form $(i(\alpha),i(\alpha))$ for
some $\alpha\in\Omega^1_X$.

The two maps $\Omega^1_X @>>> \cV_2-\cV_1$ given by $\alpha\mapsto(i(\alpha),0)$
and $\alpha\mapsto(0,-i(\alpha))$ coincide and we denote their common value
by $i$ as well. There is a short exact sequence
\[
0 @>>> \Omega^1_X @>i>> \cV_2-\cV_1 @>\pi>> \cT_X @>>> 0 \ ,
\]
where $\pi((v_2,v_1))$ is defined to be the common value of $\pi(v_2)$ and $\pi(v_1)$.

Let
\begin{multline}\label{formulas:V-minus-V}
\partial(f)=i(df),\ \ \ f(v_2,v_1)=(f*v_2,f*v_1), \\ 
[(v_2,v_1),(v'_2,v'_1)]= ([v_2,v'_2],[v_1,v'_1]),
\ \ \ \langle(v_2,v_1),(v'_2,v'_1)\rangle =
\langle v_2,v'_2\rangle+\langle v_1,v'_1\rangle
\end{multline}

It is clear that component-wise $*$-operation of $\cO_X$ is well-defined
on $\cV_2-\cV_1$. Moreover, since the associator of the $\cO_X$-action
on a vertex $\cO_X$-algebroid (i.e. the right hand side of \eqref{assoc})
does not depend on the algebroid, the operation of $\cO_X$ on $\cV_2-\cV_1$
is, in fact, associative. Since the ``anomalies'' present in \eqref{leib}
and \eqref{pairing} are algebroid-independent, $\cV_2-\cV_1$ is, in fact,
an exact Courant $\cO_X$-algebroid. It is clear that $(\cV_2-\cV_1)+\cV_1$
is canonically isomorphic to $\cV_2$.

So far we have established that, if non-empty, $\EVA_X(X)$ is a torsor
under $\ECA_X(X)$. It remains to show that $\EVA_X$ is locally
non-empty.

In the analytic or $C^\infty$ case example \ref{example:loc-pic}
provides a locally defined EVA. Indeed, locally on X there exists
an (abelian) Lie $\bbC$-subalgebra $\tau$
such that $\cT_X\cong\cO_X\otimes_\bbC\tau$.

In the algebraic setting the same example shows that $\EVA_X$ is
non-empty locally in \'etale topology. Since 
$H^2_{\text{\'et}}(X,\Omega^2_X @>d>> \Omega^{3,cl}_X)$ is canonically isomorphic
to $H^2(X,\Omega^2_X @>d>> \Omega^{3,cl}_X)$ it follows that $\EVA_X$
is non-empty Zariski-locally. 
\end{pf}

\section{Algebroids over the de Rham complex}
All of the notions of the preceeding section generalize in an
obvious way to differential graded manifolds (i.e. manifolds
whose structure sheaves are a sheaves of commutative differential graded
algebras).

For a manifold $X$ let $X^\sharp$ denote the differential graded
manifold with the underlying space $X$ and the structure sheaf $\cO_{X^\sharp}$
the de Rham complex $\DR$. In other words, $\cO_{X^\sharp} = \bigoplus_i\Omega^i_X[-i]$
(as a sheaf of graded algebras). We will denote by $\partial_{\cO_{X^\sharp}}$
the derivation given by the de Rham differential. 

\subsection{The structure of $\cT_{X^\sharp}$}
The tangent sheaf of $X^\sharp$ (of derivations of $\cO_{X^\sharp}$),
$\cT_{X^\sharp}$ is a sheaf of differential graded Lie algebras with the
differential $\partial_{\cT_{X^\sharp}} = [\partial_{\cO_{X^\sharp}},\ \ ]$
(note that $\partial_{\cO_{X^\sharp}}\in\cT_{X^\sharp}^1$).

\subsubsection{}
Let $\widetilde\cT_X$ denote the cone of the identity endomorhism
of $\cT_X$. That is, $\widetilde{\cT_X}^i=\cT_X$ for $i=-1,0$ and
zero otherwise. The only nontrivial differential is the identity map.
The complex $\widetilde\cT_X$ has the canonical structure of a
sheaf of differential graded Lie algebras.

The natural action of $\cT_X$ (respectively $\cT_X[1]$) on $\cO_{X^\sharp}$
by the Lie derivative (respectively by interior product) gives rise
to the injective map of
DGLA
\begin{equation}\label{action-tau}
\tau : \widetilde{\cT_X} @>>> \cT_{X^\sharp} \ .
\end{equation}

The action $\tau$ extends in the canonical way to a structure
of a Lie $\cO_{X^\sharp}$-algebroid on $\cO_{X^\sharp}\otimes_\bbC\widetilde\cT_X$
with the anchor map
\begin{equation}\label{tau-DR}
\tau_{\cO_{X^\sharp}} : \cO_{X^\sharp}\otimes_\bbC\widetilde\cT_X @>>> \cT_{X^\sharp}
\end{equation}
the canonical extension $\tau$ to a $\cO_{X^\sharp}$-linear map. Note that
$\tau_{\cO_{X^\sharp}}$ is surjective, i.e. the Lie $\cO_{X^\sharp}$-algebroid
$\cO_{X^\sharp}\otimes_\bbC\widetilde\cT_X$ is transitive. We denote this
algebroid by $\widetilde\cT_{X^\sharp}$.

\subsubsection{}
Let $\cT_{X^\sharp/X}\subset\cT_{X^\sharp}$ denote the normalizer of
$\cO_X\subset\cO_{X^\sharp}$. Since the action of $\cT_X[1]$ is
$\cO_X$-linear, the map $\tau$ restricts to
\[
\tau : \widetilde{\cT_X}^{-1} = \cT_X[1] @>>> \cT_{X^\sharp/X}
\]
and (the restriction of) $\tau_{\cO_{X^\sharp}}$ factors through the map
\begin{equation}\label{tau-sub}
\cO_{X^\sharp}\otimes_{\cO_X}\cT_X[1] @>>> \cT_{X^\sharp/X}
\end{equation}
which is easily seen to be an isomorphism.

Since the action $\tau$ is $\cO_X$-linear modulo
$\cT_{X^\sharp/X}$, $\tau_{\cO_{X^\sharp}}$ induces the map
\begin{equation}\label{tau-fac}
\cO_{X^\sharp}\otimes_{\cO_X}\cT_X @>>> \cT_{X^\sharp}/\cT_{X^\sharp/X}
\end{equation}
which is easily seen to be an isomorphism.

Therefore, there is an exact sequence of graded
$\cO_{X^\sharp}$-modules
\begin{equation}\label{ses-der}
0 @>>> \cO_{X^\sharp}\otimes_{\cO_X}\cT_X[1] @>>> \cT_{X^\sharp} @>>>
\cO_{X^\sharp}\otimes_{\cO_X}\cT_X @>>> 0
 \ .
\end{equation}
The composition
\[
\cO_{X^\sharp}\otimes_{\cO_X}\cT_X[1] @>>> \cT_{X^\sharp} @>{\partial_{\cT_{X^\sharp}}}>>
\cT_{X^\sharp}[1] @>>> \cO_{X^\sharp}\otimes_{\cO_X}\cT_X[1]
\]
is the identity map.

\subsection{Exact Courant $\cO_{X^\sharp}$-algebroids}

\subsubsection{Proposition}
Every exact Courant $\cO_{X^\sharp}$-algebroid admits a unique flat
connection.

\begin{pf}
Consider a Courant algebroid
\[
0 @>>> \Omega^1_{X^\sharp} @>i>> \cQ @>{\pi}>> \cT_{X^\sharp} @>>> 0
\]
Note that, since $\Omega^1_{X^\sharp}$ is concentrated in non-negative
degrees, the map $\pi^{-1} : \cQ^{-1} @>>> \cT_{X^\sharp}^{-1}$ is
an isomorphism. Since $\cT_{X^\sharp}^{-1}$ generates $\cT_{X^\sharp}$ as a
DG-module over $\cO_{X^\sharp}$, the splitting is unique if it exists.

To establish the existence it is necessary and sufficient to show
that the restriction of the anchor map to the DG-submodule of
$\cQ$ generated by $\cQ^{-1}$ is an isomorphism.

Note that the map $\tau : \widetilde\cT_X @>>> \cT_{X^\sharp}$ lifts in a
unique way to a morphism of complexes $\widetilde\tau :
\widetilde\cT_X @>>> \cQ$. It is easily seen to be a morphism of
DGLA. Let $\cQ'$ denote the $\cO_{X^\sharp}$-submodule of $\cQ$ generated by
the image of $\widetilde\tau$ (i.e. the DG $\cO_{X^\sharp}$-submodule
generated by $\cQ^{-1}$).

Since $\widetilde\tau^{-1}$ is $\cO_X$-linear it extends to the
map
\begin{equation}\label{Q-sub}
\cO_{X^\sharp}\otimes_{\cO_X}\cT_X[1] @>>> \cQ'
\end{equation}
such that the composition
\[
\cO_{X^\sharp}\otimes_{\cO_X}\cT_X[1] @>>> \cQ' @>\pi>> \cT_{X^\sharp}
\]
coincides with the composition of the isomorphism
\eqref{tau-sub}with the inclusion into $\cT_{X^\sharp}$. Therefore,
\eqref{Q-sub} is a monomorphism whose image will be denoted
$\cQ''$, and $\pi$ restricts to an isomorphism of $\cQ''$ onto
$\cT_{X^\sharp/X}$.

Since $\widetilde\tau^0$ is $\cO_X$-linear modulo $\cQ''$ it
extends to the map
\begin{equation}\label{Q-fac}
\cO_{X^\sharp}\otimes_{\cO_X}\cT_X @>>> \cQ'/\cQ''
\end{equation}
which is surjective (since $\cQ'/\cQ''$ is generated as a
$\cO_{X^\sharp}$-module by the image of $\widetilde\tau^0$), and such
that the composition
\[
\cO_{X^\sharp}\otimes_{\cO_X}\cT_X @>>> \cQ'/\cQ'' @>\pi>>
\cT_{X^\sharp}/\cT_{X^\sharp/X}
\]
coincides with the isomorphism \eqref{tau-fac}. Therefore,
\eqref{Q-fac} is an isomorphism.

Now the exact sequence \eqref{ses-der} implies that $\pi$
restricts to an isomorphism $\cQ' @>\cong>> \cT_{X^\sharp}$. The desired
splitting is the inverse isomorphism. It is obviously compatible
with brackets, hence, is a flat connection.
\end{pf}

\subsubsection{Corollary}\label{ECA-DR-final}
$\ECA_{\cO_{X^\sharp}}$ is equivalent to the final stack.

\subsubsection{Corollary}\label{corollary:existence-uniqueness}
An exact vertex $\cO_{X^\sharp}$-algebroid exists and is unique up to
canonical isomorphism.
\begin{pf}
Since $\EVA_{\cO_{X^\sharp}}$ is an affine space under $\ECA_{\cO_{X^\sharp}}$
the uniqueness (local and global) follows from Corollary
\ref{ECA-DR-final}. Local existence and uniqueness implies
global existence.
\end{pf}

\subsection{The exact vertex $\cO_{X^\sharp}$-algebroid}\label{subsection:construction}
For a vertex $\cO_{X^\sharp}$-algebroid $\cV$ with anchor map $\pi$ we will
call a $\widetilde\cT_X$-rigidification of $\cV$ a map of DG
Leibniz algebras $\psi : \widetilde\cT_X @>>> \cV$ such that the
composition $\pi\circ\psi$ coincides with the map
\eqref{action-tau}. A rigidified vertex $\cO_{X^\sharp}$-algebroid is a pair
$(\cV,\psi)$, where $\cV$ is a vertex $\cO_{X^\sharp}$-algebroid and $\psi$
is a $\widetilde\cT_X$-rigidification of $\cV$. A morphism of
rigidified vertex $\cO_{X^\sharp}$-algebroids is a morphism of vertex
$\cO_{X^\sharp}$-algebroids which commutes with respective rigidifications.

The functor which assigns to a vertex $\cO_{X^\sharp}$-algebroid $\cV$ the
set of rigidifications of $\cV$ is representable. The initial
$\widetilde\cT_X$-rigidified vertex $\cO_{X^\sharp}$-algebroid
$\widetilde\cU$ may be described as follows.
As a sheaf of $\bbC$-vector spaces
\[
\widetilde\cU =
\Omega^1_{X^\sharp}\oplus\cO_{X^\sharp}\otimes_\bbC\widetilde\cT_X \ .
\]
There is a
unique structure of a vertex $\cO_{X^\sharp}$-algebroid on $\widetilde\cU$
such that
\begin{itemize}
\item the derivation $\partial : \cO_{X^\sharp} @>>> \Omega_{\widetilde\cU}$
is given by the composition of the exterior deriviative
$d_{X^\sharp} : \DR @>>> \Omega^1_{X^\sharp}$ with the inclusion
of the first summand $\Omega^1_{X^\sharp} @>>> \widetilde\cU$,

\item the projection onto the second summand induces an
isomorphism of Lie $\cO_{X^\sharp}$-algebroids $\overline{\widetilde\cU}
@>\cong>> \DR\otimes_\bbC\widetilde\cT_X$,

\item the restriction to $1\otimes\widetilde\cT_X$ of the adjoint
action of $\widetilde\cU$ on $\Omega_{\widetilde\cU}$ coincides with the
action of $\widetilde\cT_X$ on $\Omega^1_{X^\sharp}$ by the Lie
derivative and \eqref{action-tau}.
\end{itemize}

Namely, the action of $\cO_{X^\sharp}$ is forced by \eqref{assoc}:
\[
\alpha*(\beta\otimes\xi) = \alpha\wedge\beta\otimes\xi -
\tau(\xi)(\alpha)d_{X^\sharp}\beta - \tau(\xi)(\beta)d_{X^\sharp}\alpha \ .
\]
The symmetric pairing is the unique extension of the zero pairing
on $1\otimes\widetilde\cT_X$ which satisfies \eqref{pairing}:
\[
\langle \beta_1\otimes\xi_1,\beta_2\otimes\xi_2\rangle =
-\beta_1\tau(\xi_2)(\tau(\xi_1)(\beta_2)) -
\beta_2\tau(\xi_1)(\tau(\xi_2)(\beta_1)) -
\tau(\xi_1)(\beta_2)\tau(\xi_2)(\beta_1)
\]

The Leibniz bracket on $\widetilde\cU$ is the unique extension of
the Lie bracket on $\widetilde\cT_X$ which satisfies \eqref{leib}
and \eqref{symm-bracket}.

Consider the following condition on a
$\widetilde\cT_X$-rigidification $\psi$ of a vertex $\cO_{X^\sharp}$-algebroid
$\cV$:
\subsubsection{}\label{Lin}
{\em
The component $\psi^{-1}:\cT_X @>>> \cV^{-1}$ is
$\cO_X$-linear, i.e. $\psi^{-1}(f\xi) = f*\psi^{-1}(\xi)$ for
$f\in\cO_X$ and $\xi\in\cT_X$.
}

The functor which assigns to a vertex $\cO_{X^\sharp}$-algebroid the
collection of rigidifications of $\cV$ satisfying \ref{Lin} is
represented by the appropriate quotient $\cU$ of $\widetilde\cU$.
We claim that $\cU$ is an exact vertex $\cO_{X^\sharp}$-algebroid.

Let $K$ denote the kernel of the canonical map
$\widetilde\cT_{X^\sharp}^{-1}=\cO_X\otimes_\bbC\cT_X @>>> \cT_X$.

Let
\[
\cI = \DR*K[1] + \DR*dK\subset\widetilde\cU \ .
\]
Then, $\DR*\cI=\cI$, $\langle\widetilde\cU,\cI\rangle = 0$,
and $[\widetilde\cU,\cI]\subseteq\cI$.
Moreover, $\cI$ projects isomorphically onto the kernel of \eqref{tau-DR}.

The essentially unique exact vertex $\cO_{X^\sharp}$-algebroid is $\cU = \widetilde\cU/\cI$.


\begin{thebibliography}{ABC}
\bibitem[BD]{BD} A.~Beilinson, V.~Drinfeld, Chiral algebras, preprint.
\bibitem[GMS]{GMS} V.~Gorbunov, F.~Malikov, V.~Schechtman, Gerbes of chiral
differential operators II, preprint.
\end{thebibliography}
\end{document}